\documentclass[multi]{cambridge7A}
\usepackage[UKenglish]{babel}
\usepackage[longnamesfirst,sectionbib]{natbib}
\usepackage{chapterbib}
\usepackage{amsmath,amssymb,amsthm,amsfonts}
\usepackage{enumerate,url}
\usepackage{graphicx}
\usepackage{psfrag,  epsfig}

\setcitestyle{numbers,square,comma} 

\allowdisplaybreaks[1]

\newcommand{\dW}{\tilde W} 
\newcommand{\tQ}{\breve Q}
\newcommand{\tW}{\breve W} 
\newcommand{\tM}{\breve M} 
\newcommand{\tN}{\breve N}
\newcommand{\tX}{\breve X}
\newcommand{\tU}{\breve U}
\newcommand{\tZ}{\breve Z}
\newcommand{\tE}{\breve E}
\newcommand{\tD}{\breve D}
\newcommand{\PR}{\mathbb P}
\newcommand{\E}{\mathbb E}
\newcommand{\I}{\mathcal I}
\newcommand{\J}{\mathcal J}

\providecommand*\Index[1]{#1\index{#1}}
\providecommand*\undex[1]{} 

\begin{document}
\alphafootnotes
\author[F. P. Kelly and R. J. Williams]
{F.~P. Kelly\footnotemark\
and R.~J. Williams\footnotemark}
\chapter[Heavy traffic on a controlled motorway]
{Heavy traffic on a controlled motorway}

\footnotetext[1]{Statistical Laboratory, University of Cambridge, Centre for
  Mathematical Sciences, Wilberforce Road, Cambridge CB3 0WB; 
  F.P.Kelly@statslab.cam.ac.uk}
\footnotetext[2]{Department of Mathematics, University of California at
  San Diego, 9500 Gilman Drive, La Jolla, CA 92093--0112, USA;
  williams@stochastic.ucsd.edu. Research
  supported in part by NSF Grant DMS 0906535.}
\arabicfootnotes
\contributor{Frank P. Kelly \affiliation{University of Cambridge}}
\contributor{Ruth J. Williams \affiliation{University of California at San
  Diego}}
\renewcommand\thesection{\arabic{section}}
\numberwithin{equation}{section}
\renewcommand\theequation{\thesection.\arabic{equation}}
\numberwithin{figure}{section}
\renewcommand\thefigure{\thesection.\arabic{figure}}

\begin{abstract}
Unlimited access to a motorway network can, in overloaded conditions,
cause a loss of capacity. Ramp metering (signals on slip roads to
control access to the motorway) can help avoid this loss of capacity.
The design of ramp metering strategies has several features
in common with  the design
of access control mechanisms in communication networks. 

Inspired by models and rate control mechanisms developed
for Internet congestion control,
we propose a Brownian network model as an approximate model
for a controlled motorway and consider it operating under a proportionally
fair ramp metering policy.
We present an analysis of the performance of this model.
\end{abstract}

\subparagraph{AMS subject classification (MSC2010)}90B15, 90B20, 60K30

\section{Introduction}
The study of
\index{heavy traffic|(}heavy traffic in queueing systems began
in the 1960s, with three pioneering papers by
\index{Kingman, J. F. C.!influence|(}Kingman~\cite{K1,K2,K3}. These papers, 
and the early work of
\index{Prohorov, Y. V.}Prohorov~\cite{PRO}, 
\index{Borovkov, A. A.}Borovkov~\cite{BOR1, BOR2} and
\index{Iglehart, D. L.}Iglehart~\cite{IGL}, 
concerned a single resource. 
Since then there has been significant interest in 
\index{network|(}networks of resources, with major advances
by
\index{Harrison, J. M.}Harrison and Reiman~\cite{HR},
\index{Reiman, M.I.}Reiman~\cite{REI}, 
Williams~\cite{WILQS} and
\index{Bramson, M. D.}Bramson~\cite{BRA}.
For discussions, further references
and overviews of the very extensive literature
on heavy traffic for  networks, 
Williams~\cite{WILSN},
\index{Dai, J. G.}Bramson and Dai~\cite{BD}, 
Harrison~\cite{HARW, HARB} and
\index{Whitt W@Whitt, W.}Whitt~\cite{WHITT} are recommended. 

Research in this area is motivated in part by 
the need to understand and control the behaviour
of communications, manufacturing and service networks,
and thus to improve their design and performance.
But researchers are also attracted by
the elegance of some of the mathematical constructs:
in particular, the multi-dimensional reflecting
\index{Brown, R.!Brownian motion|(}Brownian motions that often arise as limits.

A question that arises in a wide variety of
application areas concerns how flows through a network
should be controlled, so that the network
responds sensibly to varying conditions.
\index{road traffic|(}Road traffic was an  area  
of interest to early
\index{Newell, G. F.}researchers~\cite{NEW}, 
and more recently the question has
been studied in work on modelling
the
\index{Internet}Internet. In each of these cases the network studied
is part of a larger system: for example,
 drivers generate demand
and select their routes in ways that are responsive
to the delays incurred or expected, which depend on 
the controls implemented in the road network. 
It is important to address such interactions 
between the network and the larger system, 
and in particular to understand the signals,
such as delay, provided to the larger system.

Work on Internet congestion control
generally addresses the issue of
\index{fairness|(}fairness, since there
exist situations where a given scheme might 
maximise network throughput, for example, while denying 
access to some users.  In this area 
it has been possible to integrate ideas of fairness
of a control scheme with overall system optimization: 
indeed fairness of the control scheme is often
the means by which the right information and
incentives are provided to the 
larger
\index{Maulloo, A.}\index{Tan, D.}\index{Srikant, R.}system~\cite{KMT,SRI}. 

Might some of these ideas transfer to help our understanding
of the control of road traffic? 
In this paper we present a preliminary exploration of 
a particular topic:
\index{ramp metering|(}ramp metering. Unlimited access to a motorway network
can, in overloaded conditions,
cause a loss of
\index{network!capacity|(}capacity. Ramp metering (signals on slip roads to
control
access to the motorway) can help avoid this loss of capacity. 
The problem is one of access control, a common issue for
communication networks, and in this paper we describe a ramp metering
policy,
\it proportionally fair metering\rm, inspired by
rate control mechanisms developed for the Internet.

The organisation of this paper is as follows. In Section~\ref{asq}
we review early heavy traffic results for a single queue.
In Section~\ref{mic} we describe a model 
of Internet congestion control, which we use
to illustrate the simplifications and insights heavy traffic 
allows. 
In Section~\ref{bnm} we describe a
\index{Brown, R.!Brownian network}Brownian
network model, which both generalizes a model of
Section~\ref{asq} and arises as a heavy traffic limit
of the networks considered in Section~\ref{mic}.
Sections~\ref{mic} and~\ref{bnm} are  based on the recent results
\index{Kang, W. N.}\index{Lee, N. H.}of~\cite{KKLW, KW}.
These heavy traffic models  help us to understand the
behaviour of networks operating under policies
for sharing
\index{network!capacity|)}capacity fairly.

In Section~\ref{mcm}  we develop an approach to the
design of
\index{ramp metering|)}ramp metering flow rates informed by the earlier
Sections. 
For each of three examples, we present a Brownian
network model  operating
under a proportionally fair metering policy.
Our first example is a linear
network representing a road into a city centre with
several entry points; we then discuss a tree network,
and, in Section~\ref{rc}, a simple network where 
drivers have routing choices. Within the Brownian network
models we show that in each case
the delay suffered by a driver at an entry
point to the network can be expressed as a sum of
\index{dual variable}dual variables,
one for each of the resources to be used, and that under their
stationary distribution these
dual variables are independent  exponential random variables.
For the final example we
show that the interaction of proportionally
\index{fairness|)}fair metering
with choices available to arriving 
traffic has beneficial consequences for the performance
of the system.\index{road traffic|)}

John Kingman's initial insight, that heavy traffic reveals
the essential properties of queues,
generalises to networks, where heavy traffic allows 
sufficient simplification to make clear
the most important consequences of resource allocation
policies.\index{Kingman, J. F. C.!influence|)}

\section{A single queue} \label{asq}

In this Section we review heavy traffic results
for the  M/G/1 queue, to introduce ideas that
will be important later when we look at
\index{network|)}networks.

Consider a queue with a single server of
unit capacity at which 
customers arrive as a
\index{Poisson, S. D.!Poisson process}Poisson process of rate $\nu$.
Customers bring
amounts of work for 
the server which
are independent and identically
distributed with distribution $G$, 
and are independent of the arrival process. Assume
the distribution $G$ has mean $1/\mu$ and 
finite second moment, and 
that the load on the queue, $\rho = \nu /\mu$, satisfies
$\rho <1$. 


Let $W(t)$ be the
\index{workload|(}\emph{workload} in the queue
at time $t$; for a server of unit capacity this is 
the time it would take
for the server to empty the queue if no more arrivals
were to occur after time $t$. 

Kingman~\cite{K1} showed
that the stationary distribution of $(1-\rho)W$ is 
asymptotically exponentially distributed  as $\rho \rightarrow 1$.
Current approaches to heavy traffic generally proceed via
a weaker assumption that the cumulative arrival process of work 
satisfies a
\index{functional central limit theorem}functional central limit theorem,
and use this to show that 
as $\rho \rightarrow 1$, the appropriately normalized workload  process
\begin{equation} 
\hat W(t) = (1-\rho) \ W\left(\frac{t}{(1-\rho)^2}\right),   \qquad t \geq 0
\label{scaling}
\end{equation}
can be approximated by a reflecting
Brownian motion $\dW$  on $\mathbb{R}_+$. In
the interior $(0,\infty)$ of $\mathbb{R}_+$, 
 $\dW$ behaves as a Brownian
motion with drift $-1$ and variance determined by
the variance of the cumulative arrival process of work.
When $\dW$ hits zero, then the server may become idle; 
this is where
delicacy is needed. The stationary distribution of the
reflecting Brownian motion
$\dW$ is exponential, corresponding to Kingman's early
result. 

We note an important  consequence of the scalings appearing
in the definition~(\ref{scaling}), the
\index{snapshot principle}\it snapshot principle\rm. 
Because of the different scalings applied to space and time,
the workload is of order $(1-\rho)^{-1}$ while 
the workload can change significantly
only over time intervals of order $(1-\rho)^{-2}$. 
Hence the time taken to serve the
amount of work
in the queue is asymptotically negligible compared to the
time taken for the workload to change significantly~\cite{REI1, WHITT}.

Note that the workload $(W(t), t\geq 0)$  does not depend on
the
\index{queue!discipline|(}queue discipline (provided the discipline does not
allow idling when there is work to be done), although the
\index{waiting time}waiting time 
for an arriving  customer certainly does. 
\index{Kingman, J. F. C.}Kingman~\cite{K4} makes elegant use of the 
snapshot principle to compare stationary waiting time distributions 
under a range of queue disciplines. 

It will be helpful to develop in detail a simple example. 
Consider a
\index{Markov, A. A.!Markov process|(}Markov  process
in continuous time $(N(t), t \geq 0)$
with state space $\mathbb{Z}_+$ and non-diagonal 
infinitesimal  transition rates
\begin{align} \label{MM1}
q(n,n')=
\begin{cases}
\nu & \text{if }\; n'=n+1,\\
\mu  & \text{if }\; n'=n-1\text{ and }n>0,\\
0 & \text{otherwise}. 
\end{cases}
\end{align}
Let $\rho = \nu / \mu$. If $\rho <1$ then the Markov  process
$(N(t), t \geq 0)$
has stationary distribution
\begin{equation}
\mathbb{P}\{N^s=n\} 
= (1-\rho) \rho^n, \qquad  n=0,1,2 ,\ldots 
\label{geom} 
\end{equation}
(here, the superscript
$s$ signals that the random variable is associated with the stationary
distribution).
The Markov  process corresponds to an M/M/1 queue, at which 
customers arrive as a
\index{Poisson, S. D.!Poisson process}Poisson process of rate $\nu$, and
where customers bring an amount of work for the server
which is exponentially distributed with parameter $\mu$.

Next consider an M/G/1 queue with the
\index{processor sharing}processor-sharing
\index{queue!discipline|)}discipline (under the processor-sharing
discipline, while there are $n$ customers in the queue
each receives a proportion $1/n$ of the capacity of
the server). The  process $(N(t), t \geq 0)$ is
no longer 
\index{Markov, A. A.!Markov process|)}Markov, but it nonetheless has the same
stationary distribution as in (\ref{geom}). Moreover in the
stationary regime, given $N^s =n$, the amounts of work left 
to be completed on each of the $n$ customers in the queue
form a collection of $n$ independent random variables,
each with distribution function 
\begin{equation*}
G^*(x) = \mu \int_0^x (1-G(z)) \, dz , \quad x\geq 0, 
\end{equation*}
a distribution recognisable as that of the forward
\index{recurrence time}recurrence time in a stationary
\Index{renewal process} whose
inter-event time distribution is $G$. 
Thus the stationary distribution of $W$ 
is just that of the sum of $N^s$ independent random
variables each with distribution $G^*$, where $N^s$ 
has the
\index{Benes, V. E.@Bene\v s, V. E.}distribution~(\ref{geom})~\cite{BEN, KEL}.
Let $S$ be a random variable with distribution $G$.
Then we can deduce that the stationary distribution
of $(W, N)$ has the property that in probability 
$W^s/N^s \rightarrow  \mathbb{E}(S^2)/2\mathbb{E}(S)$, 
the mean of the distribution $G^*$,
 as $\rho \rightarrow 1$. 
For fixed $x$, under the stationary distribution for the queue, let $N_x^s$ be the number of customers in the
queue with a remaining work requirement of not more than $x$. Then, 
$N_x^s/N^s \rightarrow  G^*(x)$ in probability as $\rho \rightarrow 1$.
At the level of stationary distributions, this 
is an example of a property  called 
\index{state space collapse}state-space collapse: in heavy traffic the 
stochastic behaviour of the system  is essentially given by
$W$, with more detailed information about the system
(in this case, the  numbers of customers with various remaining
work requirements) not being necessary.

The amount of work arriving
at the queue over a period of time, $\tau$, has a compound
Poisson distribution, with a straightforwardly calculated 
mean  and variance of  $\rho \tau$ and $\rho \sigma^2\tau$
respectively,   where $\sigma^2 =  \mathbb{E}(S^2)/\mathbb{E}(S) $. 
An alternative
\index{Harrison, J. M.|(}approach~\cite{HARBOOK} is to directly model the
cumulative arrival process of work
as a Brownian motion $\tE = (\tE(t), t\geq 0)$ with matching mean and variance parameters: thus 
\begin{equation*}
\tE(t) = \rho t + \rho^{1/2} \sigma \tZ(t),  \qquad t \geq 0,
\end{equation*}
where $( \tZ(t), t \geq 0)$ is a 
standard  Brownian motion.
Let 
\begin{equation*}
\tX(t) = \tE(t) - t, \qquad  \qquad t \geq 0, 
\end{equation*}
a Brownian motion starting from the origin with 
drift $-(1-\rho)$  and variance $ \rho \sigma^2$.
In this approach we define the
queue's workload $\tW(t)$ at time $t$ by the system of equations 
\begin{eqnarray}
\tW(t) &=&\tW(0)+  \tX(t) + \tU(t),
\qquad   t  \geq   0,
\label{Q1} \\
\tU(t) &=& - 
\displaystyle{\inf_{0 \leq s \leq t}} \tX(s) , \qquad   t  \geq  0. 
\label{Q2} 
\end{eqnarray}
The interpretation of the model is as follows.
While $\tW$ is positive, it is driven
by the Brownian fluctuations caused by arrival
of work less the work served.
But when $\tW$ hits zero, the resource may
not be fully utilized. 
The process $\tU$ defined by equation~(\ref{Q2}) is continuous
and non-decreasing, and is the minimal such process that permits 
$\tW$, given by equation~(\ref{Q1}), to  remain non-negative. 
We interpret $\tU(t)$ as the cumulative unused capacity
up to time $t$.  Note that $\tU$ can increase only at times  when $\tW$
is at zero.

The stationary distribution of $\tW$
 is exponential
with mean $\rho \sigma^2/2(1-\rho)$~\cite{HARBOOK}\index{Harrison, J. M.|)}. 
This is the same as the distribution of $(1-\rho)^{-1} \dW^s$ where $\dW^s$ has the stationary 
distribution of the reflecting Brownian motion $\dW$ that  approximates the scaled process $\hat W$  given
by (\ref{scaling}).
Furthermore, the mean of the stationary  distribution of $\tW$ is
  the same as the mean of the exact stationary distribution of
the workload $W$,
calculated from its 
representation as the geometric
sum~(\ref{geom}) of independent random variables each with distribution $G^*$
and hence mean  $\mathbb{E}(S^2)/2\mathbb{E}(S) $.

In other words, for the M/G/1 queue, we obtain the same exponential stationary distribution 
either by (a) approximating the workload  arrival process directly by a Brownian motion without
any space or time scaling,
or  by (b) approximating the scaled workload process in (\ref{scaling}) by a
reflecting
Brownian motion,
finding the stationary distribution of the latter, and
then formally unwinding the spatial scaling to obtain a  distribution in the
original spatial units.
Furthermore, this exponential distribution has the same mean as the exact stationary distribution for the
workload in the M/G/1 queue
and provides a rather  good approximation, being of the same order of accuracy as 
the exponential approximation of the geometric distribution with the same mean.

The main point of the above discussion is that, 
in the context of this example, 
we observe that for the purposes of computing approximations to the stationary
\index{workload|)}workload, 
using a direct Brownian model
for the workload arrival process (by matching mean and variance parameters)
provides the same results as 
use of the
\index{heavy traffic!diffusion approximation}heavy traffic diffusion
approximation coupled with formal unwinding
of the spatial scaling,
and  the approximate stationary distribution that this
yields compares  remarkably well  with exact results.
We shall give another example of this kind of fortuitously good approximation
in Section \ref{bnm}. 
Chen and Yao
\index{Chen, H.}\index{Yao, D. D.}\cite{CY} have also noted remarkably good
results from using such `strong approximations' without any scaling.


\section{A model of Internet congestion} \label{mic}\index{network|(}

In this Section we describe a network generalization
of
\Index{processor sharing} that has been useful in
modelling flows through the
\index{Internet}Internet, and outline
a recent heavy traffic
\index{Kang, W. N.}\index{Lee, N. H.}approach~\cite{KKLW, KW} to its analysis.

\subsection{Fair sharing in a network} \label{mic1}

Consider a network with a finite set $\J$ of  \it resources\rm.
Let a \emph{route} $i$ be a non-empty
subset of $\J$, and write $j \in i$ to indicate that resource $j$
is used by route $i$. Let $\I$ be the set of
possible routes. Assume that both $\J$ and $\I$ are non-empty and finite,
and let $J$ and $I$ denote
the cardinality of the respective sets. 
Set $A_{ji}=1$ if $j \in i$, and $A_{ji}=0$
otherwise. This defines a  $J\times I$  matrix $A =(A_{ji},\, j\in \J,\, i \in \I)$ of zeroes and ones,
the
\index{resource-route incidence matrix}\it resource-route incidence matrix\rm.
Assume that
$A$ has rank $J$, so that it has full row rank.

Suppose that resource $j$ has
\index{network!capacity|(}capacity  $C_j>0$,
and that there are $n_i$
\index{network!connection|(}connections using route $i$. How might the
capacities $C=(C_j, j\in \J)$ be shared over 
the routes $\I$, given the numbers of 
connections $n=(n_i, i \in \I)$?
This is a question which has attracted attention in a
variety of fields, ranging from \Index{game theory}, through
\Index{economics} to \Index{political philosophy}. Here we describe a 
concept of
\index{fairness}fairness which is a natural extension of 
\index{Nash, J. F.}Nash's bargaining solution and, as such, satisfies certain
natural axioms of fairness~\cite{NASH}; the
concept has been used extensively in
the modelling of rate control
\index{algorithm}algorithms in the
\index{Internet}Internet
\index{Maulloo, A.}\index{Tan, D.}\index{Srikant, R.}\cite{KMT,SRI}.

Let $\I_+(n) = \{ i \in \I  : n_i>0 \}$.
A  \index{network!capacity allocation policy}\emph{capacity allocation policy}
$\Lambda = (\Lambda(n), n\in \mathbb{R}_+^I)$,
where $\Lambda(n) = (\Lambda_i(n), i \in \I)$, 
is called \emph{proportionally fair} if for each $n\in \mathbb{R}_+^I$, 
$\Lambda(n)$ solves 
\begin{align}
&\text{maximise} &&\sum_{i\in \I_+(n)} n_i\log \Lambda_i \quad &\ \label{pf1}\\
&\text{subject to} &&\sum_{i \in \I}A_{ji} \Lambda_i \leq C_j,
& j\in \J, \label{pf2}\\
&\text{over} &&\Lambda_i\geq 0,  & i \in \I_+(n), \label{pf3} \\
& && \Lambda_i = 0, & i \in \I \setminus \I_+(n).   \label{pf4}
\end{align}
Note that the constraint~(\ref{pf2}) captures the limited 
capacity of resource $j$, while constraint~(\ref{pf4}) requires
that no capacity be allocated to a route which has no 
connections. 

The problem~(\ref{pf1})--(\ref{pf4}) is
a straightforward
\index{convex optimization}convex optimization\break problem, with optimal
solution
\begin{equation} \label{cs1}
\Lambda_i(n) = \frac{n_i}{\sum_{j\in \J}q_j A_{ji}}, \qquad i \in \I_+(n), 
\end{equation}
where the variables $q=(q_j, j\in \J) \geq 0$ are
\index{Lagrange, J.-L.!Lagrange multiplier}Lagrange multipliers (or
\index{dual variable}\emph{dual variables}) for the constraints~(\ref{pf2}).
The solution to the optimization problem is unique and 
satisfies $\Lambda_i(n) >0$ for
$i\in \I_+(n)$ 
by the strict concavity on
$(\Lambda_i > 0, i\in \I_+(n))$
and boundary behaviour of the
objective function in (1.6)~\cite{KW}. 

The dual variables $q$ are unique if $n_i >0$ for all $i \in \I$,
but may not be unique otherwise. In any event they
satisfy the
\index{complementary slackness}\emph{complementary slackness} conditions
\begin{equation}\label{cs2}
q_j \left( C_j - \sum_{i\in \I}A_{ji}\Lambda_i(n) \right)
=0, \qquad
j\in \J.
\end{equation}

\subsection{Connection level model} \label{mic2}

The allocation $\Lambda(n)$ describes how capacities are shared
for a given number of connections $n_i$ on each route $i \in \I$.
Next we describe a stochastic
\index{Massouli\'e, L.}\index{Roberts, J.}model~\cite{MaRo98}
for how the number of connections within the network varies.

A connection on route $i$ corresponds to continuous transmission
of a document through the resources used by route $i$. Transmission
is assumed to occur simultaneously through all the resources
used by route $i$. Let the number of connections on route $i$ 
at time $t$ be denoted by $N_i(t)$, and let $N(t) = (N_i(t), i \in \I)$. 
We consider a
\index{Markov, A. A.!Markov process|(}Markov process in continuous time $(N(t), t \geq 0)$ 
with state space $\mathbb{Z}_+^I$ and non-diagonal 
infinitesimal transition rates 
\begin{align} \label{sflm}
q(n,n')=
\begin{cases}
\nu_i & \text{if }\; n'=n+e_i,\\
\mu_i \Lambda_i(n) & \text{if }\; n'=n-e_i\text{ and }n_i>0,\\
0 & \text{otherwise}, 
\end{cases}
\end{align}
where $e_i$ is the $i$-th unit vector in $\mathbb{Z}_+^I$, and
$\nu_i$, $\mu_i >0$, $i \in \I$.

The Markov process $(N(t), t \geq 0)$ 
corresponds to a model where new
\index{network!connection|)}connections arrive on
route $i$ as a
\index{Poisson, S. D.!Poisson process}Poisson process of rate $\nu_i$, and a
connection
on route $i$ transfers a document whose size is exponentially
distributed with parameter $\mu_i$. In the case where $I=J=1$ 
and $C_1=1$, 
the transition rates~(\ref{sflm}) reduce to 
the rates~(\ref{MM1}) of the M/M/1 queue. 

Define the \emph{load} on route $i$ to be $\rho_i = \nu_i / \mu_i$
for $i\in \I$.
It is
\index{Bonald, T.|(}\index{Massouli\'e, L.}\index{De Veciana, G.}\index{Lee, T. J.}\index{Konstantopoulos, T.}known~\cite{BM, DLK}
that the Markov  process is positive recurrent provided
\begin{equation} \label{capcon}
\sum_{i\in \I } A_{ji}\rho_i  < C_j, \qquad j\in \J.
\end{equation}
These are natural constraints: the load arriving
at the network for resource $j$ must be less than
the
\index{network!capacity|)}capacity of  resource $j$, for each $j \in \J$.
Let $[\rho]$ be the $I \times I$ diagonal matrix with
the entries of $\rho =(\rho_i, i \in \I)$ on its diagonal,
and define  $\nu$, $[\nu]$, $\mu$, $[\mu]$ similarly.



Each connection on route $i$ brings with it an amount of  work for 
resource $j$ which is exponentially distributed with mean
$1/\mu_i$, for $j \in i$. The
\index{Markov, A. A.!Markov process|)}Markov  process $N$ allows
us to estimate the 
\index{workload|(}workload for each resource: define
the  \emph{workload process} by  
\begin{equation}\label{req:W}
W(t ) = A \,  [\mu]^{-1} N(t ),   \qquad  t\geq 0.
\end{equation}

\subsection{Heavy traffic} \label{HT}

To approximate the workload in a heavily loaded connection-level model by that in
a Brownian network model, we view a given connection-level
model as a member of a sequence of such models approaching the heavy traffic limit.
More precisely, we
consider a sequence of connection-level models indexed by $r$ where the
network structure, defined by $A$ and $C$, does not
vary with $r$. Each member of the sequence is a stochastic system as
described in the previous section. We append a superscript of $r$ to
any process 
or parameter associated
with the $r^{th}$ system that depends on $r$. Thus, we have
processes $N^r$, $W^r$,
and parameters $\nu^r$.
We suppose $ \mu^r = \mu$ for all $r$, so that $\rho_i^r =\nu_i^r/\mu_i$,
for each $i\in \I$. We shall assume henceforth that the following
\emph{heavy traffic} condition holds: as $r\rightarrow \infty$,
\begin{equation}
 \label{req:numu}
\nu^r\rightarrow \nu\ \ \mbox{ and }  r \left(A \rho^r -
C\right)\rightarrow
- \theta\
\end{equation}
where $\nu_j >0$ and  $\theta_j >0$ for all $j\in \J$.  
Note that (\ref{req:numu}) implies 
that $\rho^r\to \rho$ as $r\to \infty$ and that $A\rho =C$.

We define
\index{fluid scaled process|(}fluid scaled processes $\bar N^r$, $\bar W^r$
as
follows. For each $r$ and $t\geq 0$, let
\begin{equation}
\bar N^r(t)\ = N^r(r t) /r, \qquad
\bar W^r(t)\ = W^r(r t) /r.\label{rfluone}
\end{equation}

What might be the limit of the sequence  $\{\bar N^r\}$ as $r\to
\infty$?  From the transition rates~(\ref{sflm})
and the observation that $\Lambda(rn) = \Lambda(n)$ for $r>0$,  
we would certainly expect that the limit satisfies 
\begin{equation}
\frac{d}{dt}\, n_i(t)=\nu_i -\mu_i \Lambda_i(n(t)),   \qquad  
i \in \I,
\label{eq:diff}
\end{equation}
whenever $n$ is differentiable at $t$ and $n_i(t)>0$ for all $i \in \I$.  
Indeed, this forms part of the fluid model developed in~\cite{KW} as a
\index{functional law of large numbers}functional-law-of-large-numbers approximation. Extra care is
needed in defining the fluid model at any time $t$ when $n_i(t)=0$, for any
$i \in \I$: 
the function $\Lambda(n)$ may not be continuous on the
boundary of the region $\mathbb{R}_+^I$, and so 
when any component $\bar N^r_i(t)$ is hitting 
zero, $\Lambda(\bar N^r(t))$ may jitter.

It is shown in~\cite{KW}
that the set of invariant states for the fluid model is
\begin{eqnarray*}
 \mathcal N =\left \{ n\in \mathbb{R}_+^I: n_i = \rho_i
\sum_{j \in \J} q_j A_{ji}, i\in  \I \
\hbox{ for some } q \in \mathbb{R}_+^J\right\}
\end{eqnarray*}
as we would expect  from formally setting the 
derivatives in~(\ref{eq:diff}) to zero and using relation~(\ref{cs1}).
Call ${\mathcal N}$ the
\index{invariant manifold}\it invariant manifold\rm. If $n \in {\mathcal N}$,
then since $A$ has full row rank 
the representation of $n$ in terms of $q$ is unique; 
furthermore, $\Lambda_i(n)=\rho_i $ for $i\in \I_+(n)$ and then since $A\rho =C$, the vector $q$ satisfies equation~(\ref{cs1})
and the 
\Index{complementary slackness} conditions~(\ref{cs2}), and hence  gives
\index{dual variable}dual variables for the
\index{convex optimization}optimization problem~(\ref{pf1})--(\ref{pf4}).

For each $n\in \mathbb{R}_+^{I}$, define $w(n)=(w_j(n), j\in \J)$,
the  workload  associated with $n$,  by
$w(n) =A [\mu]^{-1} n$.
For each $w\in \mathbb{R}_+^{J}$, define $\Delta (w)$ to be the unique
value of $ n\in \mathbb{R}_+^{I}$ that solves the following optimization
problem:
\begin{equation*}
\begin{array}{ll}
\hbox{minimize}&{F(n)}\\
\hbox{subject to}&\displaystyle{\sum_{i\in \I} 
A_{ji} \frac{n_i}{\mu_i}  \geq w_j,
\qquad  j \in \J ,} \\
     \hbox{over}& n_i \geq 0, \qquad  i\in \I,
     \end{array}
\end{equation*}
where
\begin{equation*}
F(n) = \sum_{i \in \I}  \frac{n_i^2}{\nu_i}, \qquad
n\in \mathbb{R}_+^{I}. 
\end{equation*}
The function $F(n)$ was introduced
\index{Bonald, T.|)}\index{Massouli\'e, L.}in~\cite{BM} and
can be used to show positive 
recurrence  of $N$ under conditions~(\ref{capcon}). 
In~\cite{KW} the difference $F(n) - F( \Delta (w(n)))$ is used
as  a
\index{Lyapunov, A. M.!Lyapunov function}Lyapunov function to show 
that any
\index{fluid scaled process|)}fluid model solution $(n(t), t \geq 0)$
converges towards the invariant manifold $\mathcal N$. 
It is straightforward to check that  $n \in  \mathcal N$ if and only if 
$n = \Delta (w(n))$ and 
it turns out that 
\begin{equation*}
\Delta (w) = [\rho] A'(A[\mu]^{-1} [\nu] [\mu]^{-1} A')^{-1} w.
\end{equation*}

Note that if $\bar N^r$ lives in the space
${\mathcal N}$ then $\bar W^r$, given by equations~(\ref{req:W}) 
and~(\ref{rfluone}) as $A  \, [\mu]^{-1}\bar N^r$, lives in
the space ${\mathcal W} = A \, [\mu]^{-1}{\mathcal N}$, which we can
write as 
\begin{equation}
\label{wcone}
{\mathcal W} = 
\Big\{ w \in \mathbb{R}_+^{J}:
w = A [\mu]^{-1} [\nu] [\mu]^{-1} A' q   \
\hbox{ for some } q\in \mathbb{R}_+^{J} \Big\}, 
\end{equation}
generally a space of lower dimension. 
Call ${\mathcal W}$ the
\index{workload!cone|(}\it workload cone\rm. 
Let
\begin{equation}  \label{Wsupj}
\begin{split}
\mathcal W^j &= \Big\{ w \in \mathbb{R}_+^{J}: 
w = A [\mu]^{-1} [\nu] [\mu]^{-1} A' q  \\ 
&\qquad \qquad
\hbox{for some } q\in \mathbb{R}_+^{J}
\hbox{ satisfying } q_j=0\Big\},
\end{split}
\end{equation}
which we refer to as  the $j^{th}$
face of the workload cone $\mathcal W$.

We define
\index{diffusion!diffusion scaled process}diffusion scaled processes  $\hat N^r$, $\hat W^r$
as follows. For each $r$ and $t\geq 0$, let
\begin{equation*}
\hat N^r (t) =  \frac{N^r(r^2t)}{r}, \qquad
\hat W^r(t) =  \frac{W^r (r^2t)}{r}.
\end{equation*}
In the next sub-section we outline the
\undex{weak convergence}convergence in distribution of the sequence
$\{(\hat N^r, \hat W^r)\} $ 
as $r\to \infty$. As preparation, note 
that if $N(t) \in  \mathcal N$ and
$N_i(t)>0$ for all $i \in \I$,
then $\Lambda_i(N(t))=\rho_i$  for all $i \in \I$.
Suppose, as a thought experiment, that for each $i \in \I$ the component
$\hat N^r_i$ behaves as the queue-length process in an independent M/M/1 queue, with 
a server of capacity $\rho_i$. Then a Brownian approximation to 
$\hat N^r_i$ would have variance $\nu_i + \mu_i \rho_i
= 2 \nu_i$. Next observe that if the covariance matrix of
$\hat N^r$ is $2[\nu]$ then the covariance 
matrix of $\hat W^r=A  \, [\mu]^{-1}\hat N^r$ is 
\begin{equation}
\Gamma = 2A \, [\mu]^{-1} [\nu] [\mu]^{-1}   A'.
\label{gamma}
\end{equation}


\section{A Brownian network model} \label{bnm}\index{Brown, R.!Brownian network|(}

Let $(A, \nu, \mu)$ be as in Section~\ref{mic}: thus
$A$ is a  matrix of zeroes and ones of 
dimension $J \times I$ and of full row rank, 
and $\nu$, $\mu$ are
vectors of positive entries of 
dimension $I$. Let  $\rho_i = \nu_i /\mu_i$, $i \in \I$, and
$\rho = (\rho_i, i \in \I)$. Let ${\mathcal W}$ and 
$\mathcal W^j$ be defined by 
expressions~(\ref{wcone}) and~(\ref{Wsupj}) respectively. 
Let $\theta \in \mathbb{R}^J$ and $\Gamma $ be given by (\ref{gamma}).


In the following, all processes are assumed to be defined
on a fixed filtered probability space $(\Omega, {\mathcal F}, \{\mathcal F_t\}, \PR)$ 
and to be adapted to the
\Index{filtration} $\{\mathcal F_t\}$.
Let $\eta$ be a probability distribution on $\mathcal W$.
Define a Brownian network model 
by the following relationships:
\begin{enumerate}[(i)]
\item  $\tilde W(t)= \tilde W(0) + \tilde X(t)+
\tilde U(t)$ for all $t\geq 0$,
\item $\tilde W$ has continuous paths, $\tilde W(t)\in
\mathcal W$ for all $t\geq 0$, and $\tilde W(0)$ has distribution $\eta$,
\item $\tilde X$ is a $J$-dimensional
Brownian motion starting from the origin with drift $-\theta$ and covariance matrix $\Gamma$
such that $\{\tilde X(t) +\theta t, \mathcal F_t, t\geq 0\}$ is a
\Index{martingale} under $\PR$,
\item for each $j\in \J$, $\tilde U_j$ is a
one-dimensional process such that
\begin{enumerate}[(a)]
\item
 $\tilde U_j$ is continuous and
non-decreasing, with $\tilde U_j(0)=0$,
\item $\tilde U_j(t)=\int_{0}^t 1_{\{\tilde W(s)\in
\mathcal W^j \}} \,d\tilde U_j(s)$ for all $t\geq 0$.
\end{enumerate}
\end{enumerate}

The interpretation of the above Brownian network model is as 
follows. In the interior of the workload cone $\mathcal W$
each of the $J$ resources are fully utilized, route $i$ 
is receiving a
\index{network!capacity}capacity  allocation $\rho_i$ for each
$i \in \I$, and the workloads $\tilde W$ are driven by 
the Brownian  fluctuations caused by arrivals
and departures of
\index{network!connection|(}connections. But when $\tilde W$ hits 
the $j^{th}$
face of the workload cone $\mathcal W$, 
resource  $j$ may not be fully utilized.
The cumulative unused capacity $\tilde U_j$ at
resource $j$ is non-decreasing, and can  increase 
only on the $j^{{\rm th}}$ face of the
\index{workload!cone|)}workload cone $\mathcal W$.

The work of Dai and
\index{Dai, J. G.}Williams~\cite{DW} establishes the 
existence and uniqueness in law of the above
\index{diffusion!diffusion process}diffusion 
$\tilde W = (\tilde W(t), t\geq 0)$. 
\index{Kang, W. N.}\index{Lee, N. H.}In~\cite{KKLW} it is shown that, if
$\theta_j >0$ 
for all $j \in \J$, then
$\tilde W$  has a unique stationary distribution; furthermore,
if $\tilde W^s$ denotes a random variable with this stationary distribution, then the 
components of $\tilde Q^s = 2 \Gamma^{-1}  \tilde W^s$ 
are independent and $\tilde Q^s_j$ is exponentially distributed with
parameter $\theta_j$ for each $j \in \J$.  

Now let $C$ be
a vector of positive entries of dimension $J$, define a
sequence of networks as in Section~\ref{HT}, and suppose  
$\theta$ and $C$  are related by the  heavy traffic
condition~(\ref{req:numu}). 
In~\cite{KKLW} it is shown that, subject to 
a certain local traffic condition on the matrix $A$
and suitable convergence of
initial variables $(\hat W^r (0), \hat N^r(0))$, 
the pair $(\hat W^r,\hat N^r)$ converges in distribution as $r\rightarrow
\infty$ to a continuous
process $(\tilde W, \tilde N)$ where $\tilde W$ is the
above diffusion and $\tilde N=\Delta (\tilde W)$.
The proof in~\cite{KKLW} relies on both the existence
and uniqueness results of~\cite{DW}
and an associated
invariance principle developed by Kang and Williams~\cite{KaWi}.
(The local traffic condition under which convergence is 
established requires that the matrix $A$ contains amongst 
its columns the columns of the $J \times J$ identity matrix:
this corresponds to each resource serving at least one
route which uses only that resource. The local traffic
condition is not needed to show that
$\tilde W$ has the aforementioned stationary
distribution; that requires only the weaker condition
that $A$ have full row rank.)

It is convenient to define $\tilde Q = 2 \Gamma^{-1}  \tilde W$, a
process of
\index{dual variable}\it dual variables. \rm 
From this, the form of $\Delta $, and the relation $\tilde N=\Delta (\tilde W)$,
it follows that
$\tilde N = [\rho]  A' \tilde Q$. 
The dimension of 
the space in which $\tilde Q$ lives is $J$, 
and so this is an example of
\index{state space collapse}state-space collapse, with the $I$-dimensional process $\tilde N$
living on a $J$-dimensional
manifold where $J\leq I$ is often considerably less than $I$. 

Using the stationary distribution for $\tilde W$, we see that  
$\tilde N^s =[\rho]A'\tilde Q^s $ has the stationary distribution
of $\tilde N$.
Then, after formally 
unwinding the spatial scaling used to obtain our Brownian approximation, we obtain
the following simple approximation
for the stationary distribution of the number-of-connections
process in the  original model described in
Section \ref{mic2}:
\begin{equation} \label{approx}
N^s_i \approx \rho_i \sum_{j \in  \J} Q^s_j A_{ji}, \quad i\in \I,
\end{equation}
where $Q^s_j$, $j \in  \J$, are independent and $Q^s_j$ is
exponentially distributed with parameter
$C_j - \sum_{i\in \I} A_{ji} \rho_i$.


As mentioned  in Section~\ref{asq},  an alternative approach is to directly model
the cumulative arrival process of work  for each route $i$ as a Brownian motion:
\begin{equation*}
\tE_i(t) = \rho_i t +\left(\frac{2\rho_i}{\mu_i}\right)^{1/2} \tZ_i(t) , \quad t\geq 0,
\end{equation*}
where $(\tZ_i(t), t\geq 0)$, $i\in \I$, are independent standard 
Brownian motions; here  the form of the variance parameter takes account of the fact that
the document sizes are exponentially distributed.
Under this model, the potential  \Index{netflow}  (inflow minus potential outflow,   ignoring
underutilization of resources) process of work
for resource $j$ is
\begin{equation*}
\tX_j (t) = \sum_{i\in \I} A_{ji} \tE_i(t) - C_j t , \quad t\geq 0, 
\end{equation*}
a $J$-dimensional
Brownian motion starting from the origin with drift $A\rho-C$
and covariance matrix 
$2 A [\rho] [\mu]^{-1}A'=\Gamma$.
Then the   workload is modelled  by 
a $J$-dimensional process $\tW$ that satisfies properties (i)--(iv) above,
but with $\tW $ in place of $\tilde W$ and
$A\rho -C$ in place of the drift $-\theta $; the covariance matrix remains the same.
By the results of
\index{Kang, W. N.}\index{Lee, N. H.}\cite{KKLW}, if $A\rho <C$,  there is a
unique stationary distribution for the process
$\tW$ such that if  $\tW^s$ has  this stationary distribution then 
the components of $\tQ^s =2\Gamma^{-1} \tW^s$ are independent  and $\tQ^s_j$ is exponentially distributed
with parameter $C_j-\sum_{i\in \I} A_{ji}\rho_i$ for each $j\in \J$.
The random variable
\begin{equation}\label{approxp}
\tN^s= [\rho] A'\tQ^s ,
\end{equation}
has the stationary distribution of 
$\tN =[\rho]A' \tQ$, which is the same as 
the distribution of the right member of (\ref{approx}).
Thus,  just as in the simple case considered in Section \ref{asq}, 
in this connection-level model, 
using the direct Brownian model
yields the same approximation for the stationary distribution
of the
\index{network!connection|)}number-of-connections process as that obtained 
using the
\index{heavy traffic!diffusion approximation}heavy traffic diffusion
approximation and formally unwinding
the spatial scaling in its stationary distribution.

If we specialize the direct Brownian network model 
to the case where $I=J=1$ and $C=1$, then
we obtain the Brownian model of Section~\ref{asq}, with
$\Gamma = 2 \nu/\mu^2 = \rho \sigma^2$ and where
the stationary distribution for $\tW$ is exponentially
distributed with mean $\rho \sigma^2/2(1-\rho)$,
yielding the same approximation as in
 Section \ref{asq}.

A more interesting example is obtained when $I=J+1$
and $A$ is the $J\times (J+1)$ matrix: 
\begin{equation*}
 A = \left( \begin{array}{ccccc}
1 &  0  & \ldots & 0&1 \\
0  &  1 & \ldots& 0&1\\
\vdots & \vdots & \ddots & \vdots&\vdots\\
0&0&\ldots  & 1&1
\end{array} \right) , 
\end{equation*}
so that 
$J$ routes 
each use a single resource in such a way that there is
exactly one such 
route for each resource, and 
one route uses all $J$ resources.
In this case, the stationary distribution
given by (\ref{approxp})  
accords remarkably well with the exact
stationary distribution  described by
\index{Massouli\'e, L.}\index{Roberts, J.}Massouli\'e and Roberts
\cite{MaRo98}; it is
again of the order of accuracy of the exponential approximation of the geometric
distribution with the same mean. (We refer the interested reader to
\index{Kang, W. N.}\index{Lee, N. H.}\cite{KKLW} for the details
of this good approximation.)


In this Section and in Section \ref{asq}
we have seen intriguing examples of remarkably 
good approximations
that the direct Brownian modelling approach 
can yield. 
Inspired by this,
in the next two Sections we explore the use of the 
direct Brownian network model 
as a representation of workload
for a controlled motorway.  
Rigorous
justification for use of 
this modelling framework in the motorway context has yet to be
investigated. See the last section of the
paper  for further comments on this issue.

\section{A model of a controlled motorway} \label{mcm}\index{road traffic|(}

Once motorway traffic exceeds a certain threshold level (measured in terms
of density---the number of vehicles per mile) both vehicle speed and 
vehicle throughput drop
\index{Gibbens, R. J.}\index{Werft, W.}\index{Saatci, Y.}\index{Varaiya, P.}precipitously~\cite{GW, GWa, VAR1}.
The smooth pattern of flow that existed at lower densities breaks down, and
the driver experiences stop-go traffic. Maximum vehicle throughput (measured
in terms of the number of vehicles per minute) occurs at quite high
speeds---about 60 miles per hour on Californian freeways and on London's
orbital motorway, the M25~\cite{GW, GWa, VAR1}---while after flow breakdown
the average speed may drop to 20--30 miles per hour. Particularly problematic
is that flow breakdown may persist long after the conditions that provoked its
onset have disappeared.

Variable speed limits lessen the number and severity of accidents on
congested roads and are in use, for example,
 on the south-west quadrant of the M25.
But variable speed limits
do not avoid the loss of throughput caused by too high
a density of
\index{Abou-Rahme, N.}\index{Harbord, B.}\index{White, J.}vehicles~\cite{ABO, HAR}.
\index{ramp metering|(}Ramp metering (signals on slip roads to
control access to the motorway) can limit the density of
vehicles, and thus can avoid the loss of
\index{Levinson, D.}\index{Zhang, L.}\index{Papageorgiou, M.}\index{Kotsialis,
A.}throughput~\cite{LZ, PAP, VAR2, ZL}.
But a cost of this is queueing delay on the approaches to the motorway.
How should ramp  metering flow rates be chosen to control
these queues, and to distribute queueing delay
\index{fairness}fairly
over the various users of the motorway?  
In this  Section we introduce a modelling approach 
to address this question, based on several of the simplifications
that we have seen arise in heavy traffic.

\subsection{A linear network}  \label{aln}

\begin{figure}[!t] 
\begin{center}
\psfrag{W1}[t][t]{$w_1$}
\psfrag{W2}[t][t]{$w_2$}
\psfrag{0}[t][t]{$0$}
\psfrag{User}[t][t]{users}
\psfrag{long}[t][t]{long, uncongested link}
\psfrag{short}[t][t]{short,}
\psfrag{C1}[t][t]{$C_1$}
\psfrag{C2}[t][t]{$C_2$}
\psfrag{C3}[t][t]{$C_3$}
\psfrag{C4}[t][t]{$C_4$}
\psfrag{x1}[t][t]{$m_1$}
\psfrag{x2}[t][t]{$m_2$}
\psfrag{x3}[t][t]{$m_3$}
\psfrag{x4}[t][t]{$m_4$}

\psfrag{short1}[t][t]{congested links}
\psfrag{Cache1?}[t][t]{cache?}
\psfrag{Cache2?}[t][t]{cache?}
\centerline\leavevmode\epsfig{file=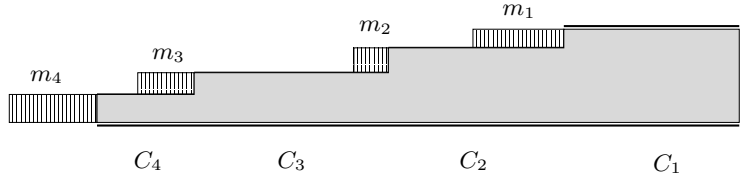,width=3.8in}

\end{center}
\caption{
Lines of size $m_1$, $m_2$, $m_3$, $m_4$ are held on the slip roads leading
to the main carriageway. Traffic on the main carriageway is free-flowing.
Access to the main carriageway from the slip roads is metered, so that 
the capacities $C_1$, $C_2$, $C_3$, $C_4$ of successive sections are not overloaded.
}
 \label{fig:road}
\end{figure}

Consider the linear\footnote{We caution the reader 
that here
we use the descriptive term `linear network' in a manner that differs from
its use in 
 \cite{KKLW}.} road network
illustrated in Figure~\ref{fig:road}.
Traffic can enter the main carriageway from lines at entry points, and
then travels from left to right, with all traffic destined for
the exit at the right hand end (think of this as a model
of a road collecting traffic all bound for a city).
Let $M_1(t)$, $M_2(t)$, \ldots, $M_J(t) $ 
taking values in $\mathbb R_+$
be the line sizes\footnote{The term
line size  is used here 
to mean a quantity measuring the
amount of work in the queue, rather than
the more restrictive number of jobs that is often
associated with the term queue size.} 
at the entry points at time $t$, and
let $C_1$, $C_2$, \ldots, $C_J$ be the respective capacities of sections
of the road. We assume the road starts at the left hand end, with
line $J$ feeding an initial section of
\index{network!capacity|(}capacity $C_J$, and
that $C_1> C_2>\ldots >C_J>0$.
The corresponding
\Index{resource-route incidence matrix} is the square matrix 
\begin{equation}
 A = \left( \begin{array}{cccc}
1 &  1  & \ldots & 1 \\
0  &  1 & \ldots & 1\\
\vdots & \vdots & \ddots & \vdots\\
0  &   0       &\ldots & 1
\end{array} \right) . 
\label{Asquare}
\end{equation}

We model the traffic, or work, 
arriving at line $i$, $i \in \I$,  as follows: let
$E_i(t)$ be the
cumulative inflow to line $i$ over the time interval $(0,t]$, and
assume $(E_i(t), t\geq 0)$ is an ergodic process with 
non-negative, stationary increments, with  $\E[E_i(t)] = \rho_i t$, 
where $\rho_i>0$,
and suppose
these processes are independent over $i \in \I$. 
Suppose
the metering rates for 
lines 1, 2, \ldots , $J$ at time $t$ can be chosen to be any
measurable vector-valued function $\Lambda = \Lambda(M(t))$ 
satisfying constraints~(\ref{pf2})--(\ref{pf4}) with 
$n=M(t)$, and such that 
\begin{equation} 
M_i(t) = M_i(0) +E_i(t) -  \int_0^t \Lambda _i (M(s)) \, ds   \geq 0  ,
\qquad t \geq 0  
\label{eeeT}
\end{equation} 
for $i \in \I$. 
Observe that we do not take into account travel time along the
road: motivated by the 
\index{snapshot principle}snapshot principle, we suppose 
that $M(\cdot)$ varies relatively slowly compared with the time
taken to travel through the system.\footnote{The time
taken for a vehicle to travel through the system comprises
both the queueing time at the entry point and the travel
time along the  motorway. If the motorway is free-flowing, 
the aim of ramp metering, then the  travel
time along the  motorway may be reasonably modelled by a
constant not dependent on $M(t)$, say $\tau_i$ from entry point $i$.
A more refined treatment 
might insist that the rates $\Lambda(M_i(t-\tau_i), i \in \I)$ satisfy
the capacity constraints~(\ref{pf2}).
We adopt the
simpler approach, since we expect that in
heavy traffic travel times along the motorway 
will be small compared with the time taken for 
$M(t)$ to change significantly.}



How might the rate function $\Lambda(\cdot )$ be chosen? We begin by a discussion of
two extreme strategies. First we consider a strategy that prioritises
the upstream entry points. 
Suppose the metered rate from  line $J$, $(\Lambda_J(M(t)), t \geq 0)$, is
chosen so that for each $t \geq 0$ the cumulative outflow from 
line $J$, $\int_0^t \Lambda_J(M(s)) \, ds$, is 
maximal, subject to the constraint~(\ref{eeeT}) and 
$\Lambda_J(M(t)) \leq  C_J$ for all $t\geq 0$: 
 thus
there is equality in the latter constraint whenever $M_J(t)$ is positive.
For each of $j=J-1$, $J-2$, \ldots, 1 in turn 
define $ \int_0^t\Lambda_j(M(s))\, ds$ to be maximal,
subject to the constraint~(\ref{eeeT}) and 
\begin{equation} 
\Lambda_j(M(t)) \leq   C_j - \sum_{i=j+1}^J \Lambda_i(M(t)),\qquad t \geq 0.
\label{eeT}
\end{equation}
In consequence there is equality in constraint~(\ref{eeT}) at time $t$ if
$M_j(t)>0$, and by induction for each $t \geq 0$ the cumulative flow along 
link $j$,
$\int_0^t \sum_{i=j}^J \Lambda_i(M(s))  \, ds$, 
is maximal, for $j=J$, $J-1$, \ldots, 1. Thus this strategy minimizes, for all times $t$,
the sum of the line sizes at time $t$, $\sum_{j} M_j(t)$. 


The above optimality property is compelling if the arrival
patterns of traffic are exogenously determined.
The strategy will, however, concentrate  delay upon the 
flows entering the system at the more downstream
entry points. This seems intuitively
\index{fairness|(}unfair, since these
flows use fewer of the system's resources, and it 
may well have perverse and suboptimal consequences if
it encourages growth in the load $\rho_i$ arriving at 
the upstream entry points. For example, growth in  $\rho_J$ 
 may cause the natural constraint~(\ref{capcon}) to be violated,
even while traffic arriving at line $J$ suffers only 
a small amount of additional delay.

Next we consider a strategy that  prioritises
the downstream entry points. To present the argument most
straightforwardly, let us suppose that 
the cumulative inflow to  line $i$
is discrete, i.e., $(E_i(t), t \geq 0)$ is constant except at an
increasing, countable sequence of times $t \in (0, \infty)$,
for each $i \in \I$.
Suppose the inflow from line $1$ is chosen
to be $\Lambda_1(M(t)) = C_1$ whenever $M_1(t)$ is positive, and zero 
otherwise. Then link $1$ will be fully utilized by the inflow
from line $1$ a proportion $\rho_1/C_1$ of the time. Let
$\Lambda_j(M(t)) = C_j$ whenever both 
$M_j(t)$ is positive and $\Lambda_i(M(t)) = 0$ for 
$i<j$, and let $\Lambda_j(M(t)) = 0$ otherwise. 
This strategy minimizes lexicographically the vector $(M_1(t), 
M_2(t), \ldots , M_J(t))$ at all times $t$. 
Provided the system
is stable, link $1$  will be utilized solely by the  inflow
from line $j$ a proportion $\rho_j/C_j$ of the time. 
Hence the system
will be unstable if 
$$
\sum_{j=1}^J \frac{\rho_j}{C_j} > 1, 
$$
and thus may well be unstable even when the
condition~(\ref{capcon}) is satisfied. Essentially the strategy
starves the downstream links, preventing them from
working at their full
capacity. Our assumption that the cumulative
inflow to line $i$ is discrete is not essential for this
argument: the stability region will be reduced from~(\ref{capcon})
under fairly general conditions. 

The two extreme strategies we have described each have their
own interest: the first has a certain optimality property
but distributes delay unfairly, 
while the second can destabilise a network even when
all the natural capacity constraints~(\ref{capcon}) are
satisfied. 


\subsection{Fair sharing of the linear network}

In this sub-section we describe our preferred ramp metering policy
for the linear network, 
and our Brownian network model for its performance.

Given the line sizes 
$M(t)=m$, we suppose the metered rates $\Lambda(m)$ are
chosen to be proportionally fair: that is, the
\index{network!capacity allocation policy}capacity allocation policy
$\Lambda(\cdot)$  solves the
\index{convex optimization}optimization problem~(\ref{pf1})--(\ref{pf4}). 
Hence for the linear network we have from relations~(\ref{cs1})--(\ref{cs2})
that 
\begin{equation*}
\Lambda_i(m) =  \frac{m_i}{\sum_{j=1}^{i} q_j }, \qquad i \in \I_+(m), 
\end{equation*}
where the $q_j$ are
\index{Lagrange, J.-L.!Lagrange multiplier}Lagrange multipliers satisfying
\begin{equation}
q_j\geq 0,\quad q_j \left( C_j - \sum_{i=j}^J \Lambda_i(m) \right)
=0, \qquad
j\in \J.
\label{cs2a}
\end{equation}
Under this policy 
the total flow along section $j$ will be its
\index{network!capacity|)}capacity $C_j$
whenever $q_j>0$. 

Given line sizes $M(t)=m$, the ratio $m_i/ \Lambda_i(m)$ is the
time it would take
to process the work currently
in  line $i$ at the current metered rate for line $i$. Thus
\begin{equation}
d_i= \sum_{j=1}^{i} q_j ,   \qquad i \in \I,
\label{delay}
\end{equation}
give estimates, based on current line  sizes, of queueing delay in each
of the $I$ lines. 
Note that
these estimates do not take into account any change 
in the line sizes over the time taken for work to 
move through the line. 

Next we describe our direct Brownian network model for the linear network
operating under the above policy.  
We make the assumption
that the inflow to line $i$ is  a
Brownian motion $\tE_i=(\tE_i(t),t\geq 0)$ starting from
the origin with drift  $\rho_i$ and variance parameter $\rho_i \sigma^2$,
and so can be written in the form 
\begin{equation}
\tE_i(t) = \rho_i t + \rho_i^{1/2} \sigma \tZ_i(t)
 , \qquad t \geq 0 
\label{eT}
\end{equation}
for $i \in \I$, where $(\tZ_i(t), t \geq 0)$, $i \in \I$, are independent
 standard Brownian motions. 
For example, if the inflow to each line were a
\index{Poisson, S. D.!Poisson process}Poisson process, then 
this would be the central limit approximation, with $\sigma =1$. 
More general choices of $\sigma$ could arise from either
a compound Poisson process, or the central limit
approximation to a large class of inflow processes. 

Our Brownian network model will be a generalization of the
model (\ref{Q1})--(\ref{Q2}) of a single queue, 
and a specialization of the  
model of Section~\ref{bnm} to the case 
where $\mu_i = 2/\sigma^2$, $i \in \I$, and the matrix $A$ is of
the form~(\ref{Asquare}). 


Let
\begin{equation*} 
 \tX_j(t) = \sum_{i \in \I}A_{ji} \tE_i(t) - C_j t, \quad t\geq 0;
\end{equation*}
note that the first term is the cumulative workload 
entering the system for resource $j$ over the interval $(0,t]$.
Write $\tX(t) = (\tX_j(t), j \in \J)$ and 
$\tX = (\tX(t), t \geq 0)$. Then
$\tX$ is a  $J$-dimensional Brownian motion starting from
the origin with drift $  A\rho-C$ and covariance 
matrix $\Gamma= \sigma^2 A \,  [\rho] A'$. 
We  assume the
stability condition~(\ref{capcon}) is satisfied, so that
$A\rho <C $.

Write
\begin{equation}
{\mathcal W} =   A [\rho] A'  \,  \mathbb{R}_+^J
\label{wc}
\end{equation}
for the
\index{workload!cone|(}workload cone, and
\begin{equation}
\mathcal W^j = \{  A \, [\rho] A'  q : q \in \mathbb{R}_+^J, q_j =0  \},
\label{wcj}
\end{equation}
for the $j^{{\rm th}}$ face of $\mathcal W$. 
Our Brownian network model for the
resource level workload $AM$ 
is then the process $\tW$ 
defined by properties (i)--(iv) of Section~\ref{bnm} with 
$\tW $ in place of $\tilde W$,  $C-A\rho$ in place of $\theta$ and 
$\Gamma = \sigma^2 A[\rho]A'$. 

The form~(\ref{Asquare}) of the matrix $A$ allows us
to rewrite the workload cone~(\ref{wc}) as
\begin{equation*}
{\mathcal W} = \left\{ w \in \mathbb{R}^J : 
 \frac {w_{j-1}-w_j}{\rho_{j-1}} \leq \frac {w_{j}-w_{j+1}}{\rho_{j}}, \,
j=1,2,\ldots ,J
\right\}
\end{equation*}
where $w_{J+1}=0$ and  we interpret the left hand side of the inequality 
as $0$ when $j=1$.
Under this model, 
at any time $t$ when the workloads $\tW(t)$ are in
the interior of the workload cone $\mathcal W$, 
each resource is fully utilized. But when
$\tW$ hits the  $j^{{\rm th}}$
face of the workload cone $\mathcal W$, 
resource  $j$ may not be fully utilized.
Our model
corresponds to the assumption that there is no more loss of
utilization 
than is necessary to prevent $\tW$ from
leaving $\mathcal W$. This assumption is made for our
Brownian network  model by analogy with the 
results reviewed in Sections~\ref{mic} and~\ref{bnm},
where it emerged as a property of the
\index{heavy traffic!diffusion approximation}heavy traffic
diffusion approximation.

In a similar manner to that in
Section \ref{bnm}, we define a process of 
\index{dual variable}dual variables: $\tQ = (A[\rho]A')^{-1} \tW$.
Since $\tW$ is our model
for $AM$, our Brownian model for the line sizes is given by
\begin{equation}
\tM = A^{-1} \tW = [\rho]  A' \tQ. 
\label{delay1}
\end{equation}
Within our Brownian model we represent 
(nominal)
delays $\tD =(\tD_i, i \in  \I)$ at
each line as given by
\begin{equation}
\tD =  A' \tQ, 
\label{delay2}
\end{equation}
since these would be the delays if line sizes remained constant 
over the time taken for a 
unit of traffic to move through the line, with $\rho_i$ 
both the arrival rate and metered rate at line $i$.\footnote{The 
nominal delay $\tD_i(t)$ for line $i$ at time $t$ 
will not in general be the
realized delay (the time
taken for the amount of work $\tM_i(t)$ found in line $i$ at
time $t$ to be metered from line $i$). 
Since  $  A\rho < C$ the
metered rate
$\Lambda_i(m)$ will in general differ from $\rho_i$  even when 
$m = A^{-1} \tW$ and $\tW \in \mathcal W$. 
Our definition of nominal delay is informed by
our earlier heavy traffic results: as $A\rho$ approaches $C$ we 
expect
scaled realized delay to converge to scaled nominal delay.
Metered rates do
fluctuate as a unit of traffic moves
through the line, but we expect 
less and less so as the system moves into heavy
traffic.
 }
Relation~(\ref{delay2}) becomes, for the linear network, 
\begin{equation*}
\tD_i = \sum_{j=1}^{i} \tQ_j ,   \qquad i=1,2,\ldots ,J, 
\end{equation*}
parallelling relation~(\ref{delay}). 
Note that when $\tW$ hits the  $j^{{\rm th}}$
face of the workload cone $\mathcal W$, then $\tQ_j=0$ 
 and 
$\tD_{j-1} = \tD_{j}$; thus the loss
of utilization at resource $j$ 
when $\tW$ hits the  $j^{{\rm th}}$
face of the workload cone $\mathcal W$
is just sufficient 
to prevent the delay at line $j$ becoming smaller than
the delay at the downstream line $j-1$.

If $\tW^s$ has the
stationary distribution of $\tW$, then
the components
of
$\tQ^s =(A[\rho]A')^{-1} \tW^s$  
are independent and $\tQ_j^s$ 
is exponentially distributed with
parameter 
\begin{equation*}
\frac{2}{\sigma^2}\left(C_j -  \sum_{i=j}^J \rho_i \right),  \qquad j=1,2,\ldots,J.
\end{equation*}
The stationary distributions of $\tM$ and $\tD$ are
then given by the distributions of $\tM^s$ and $ \tD^s$, respectively, where
\begin{equation*} 
\tM_i^s =  \rho_i \sum_{j=1}^{i} \tQ_j^s,
\quad
 \tD_i^s= \sum_{j=1}^i \tQ^s_j   \qquad i=1,2,\ldots ,J.
\end{equation*}





In the above example the matrix $A$ is invertible. As an
example of a network with a non-invertible $A$ matrix,
suppose that in the linear network 
illustrated in Figure~\ref{fig:road} one section of road
is unconstrained, say $C_3 = \infty$. Then, removing the
corresponding row from the
\Index{resource-route incidence matrix} we have 
\[ A = \left( \begin{array}{cccc}
1 &  1  & 1 & 1 \\
0  &  1 & 1 & 1\\
0  &   0       & 0 & 1
\end{array} \right).  \]
The workload cone is 
the collapse of $\mathcal W$ obtained by setting
$w_2=w_3$,
 and in consequence 
the construction of $\tW$ and $\tM=[\rho]A'(A[\rho]A')^{-1}\tW$ 
enforces the relationship $\tM_2/\tM_3 = \rho_2/\rho_3$. Since
the matrix $A$ is not invertible, this is no longer a 
necessary consequence of the network topology, but is a
natural modelling assumption,  motivated by the forms of
\index{state space collapse}state-space collapse we have seen earlier.  Essentially lines $2$ and
$3$ use the same network resources and face the same
queueing delays.

A Brownian network model of the first strategy from Section~\ref{aln}
could also be constructed, but the workload cone
and its faces 
would not be of the required form~(\ref{wc}) 
and~(\ref{wcj}), but instead
would be defined by 
\begin{equation}
{\mathcal W} = \{ w \in \mathbb{R}^J :
  0 \leq w_J \leq \ldots  \leq w_2 \leq w_1 \}, 
\label{Wineq} 
\end{equation}
and the requirement that 
if $w \in \mathcal W^j$ then $w_j=w_{j+1}$, with the
interpretation $w_{J+1}=0$. Thus face $j$ represents
the requirement that the workload for resource $j$ comprises
at least the workload for resource $j+1$, for $j=1$, 2, \ldots, $J$.
Under this model, resource $j$ is
fully utilized except when
$\tW$ hits the  $j^{{\rm th}}$
face of the workload cone~(\ref{Wineq}):
it is not possible
for $\tW$ to leave $\mathcal W$, since the
constraints expressed in the form~(\ref{Wineq}) follow
necessarily from the topology of the network embodied in $A$.
The model
corresponds to the assumption that there is no more loss of
utilization
than is a necessary consequence of the network
topology.
Note that the proportionally fair
policy may fail to fully utilize a resource not only 
when this is a necessary consequence of the network
topology, but also when this would cause an upstream
entry point to obtain more than what the policy considers
a fair share of a scarce downstream resource.



\begin{figure}[!t]
\begin{center}
\psfrag{C1}[t][t]{$C_1$}
\psfrag{C2}[t][t]{$C_2$}
\psfrag{C3}[t][t]{$C_3$}
\psfrag{C4}[t][t]{$C_4$}
\psfrag{C5}[t][t]{$C_5$}
\psfrag{C6}[t][t]{$C_6$}

\psfrag{x1}[t][t]{$m_1$}
\psfrag{x2}[t][t]{$m_2$}
\psfrag{x3}[t][t]{$m_3$}
\psfrag{x4}[t][t]{$m_4$}
\psfrag{x5}[t][t]{$m_5$}
\psfrag{x6}[t][t]{$m_6$}

\centerline\leavevmode\epsfig{file=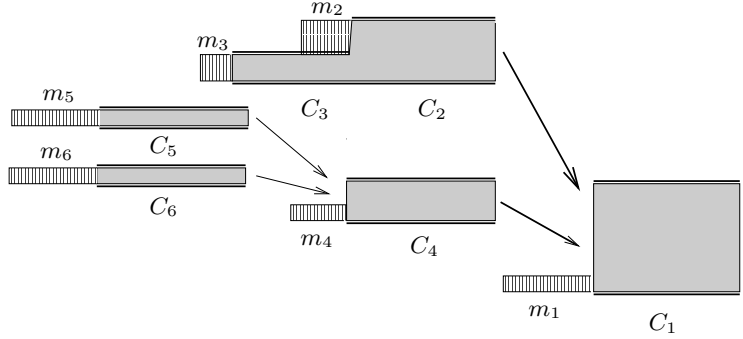,width=3.8in}

\end{center}
\caption{
There are six sources of traffic, starting in various lines and all
destined to eventually traverse section $1$ of the road. Once traffic
has passed through the queue at its entry point, it does not queue again.
}
 \label{fig:roadtree}
\end{figure}

\subsection{A tree network}

Next consider the tree network illustrated in
Figure~\ref{fig:roadtree}. Access is meter\-ed at the six entry points so 
that the capacities $C_1$, $C_2$, \ldots, $C_6$ are not overloaded. There is 
no queueing after the entry point, and the capacities 
satisfy the conditions $C_3<C_2$, $C_5+C_6<C_4$, $C_2+C_4<C_1$.

Given the  line sizes 
$M(t)=m$, we suppose the metered rates $\Lambda(m)$ are
chosen to be proportionally fair: that is, the
\index{network!capacity allocation policy}capacity allocation policy
$\Lambda(\cdot)$ solves the
\index{convex optimization}optimization problem~(\ref{pf1})--(\ref{pf4})
where for this network
\begin{equation*}
 A = \left( \begin{array}{cccccc}
1 &  1         & 1 & 1  & 1 & 1\\
0  &  1        & 1  & 0 & 0 & 0\\
0 & 0 & 1      & 0  & 0 & 0\\
0  &   0       & 0 & 1 & 1 & 1\\
0  &   0       & 0 & 0 & 1 & 0\\
0  &   0       & 0 & 0 & 0 & 1 
\end{array} \right) . 
\end{equation*}

We assume, as in the last Section, 
that the cumulative inflow
of work to  line $i$ is given by equation~(\ref{eT}) for $i \in \I$, where
$(\tZ_i(t), t \geq 0)$, $i \in \I$, are independent
 standard
Brownian motions. Our Brownian 
network model is again the process $\tW$ 
defined by properties (i)--(iv) of Section~\ref{bnm}
with $\tW$ in place of $\tilde W$, 
$C-A\rho$ in place of $\theta$,  $\Gamma =\sigma^2
A[\rho]A'$,
and with the workload cone and its faces  defined by
equations~(\ref{wc})--(\ref{wcj}) for the
above choice of $A$. We  assume the 
stability condition~(\ref{capcon}) is satisfied, so that all 
components of $ C-A\rho$ are positive.

If $\tW^s$ denotes a random variable with the
stationary distribution of $\tW$, then
the components
of
$\tQ^s =(A[\rho]A')^{-1} \tW^s$  
are independent and $\tQ_j^s$ 
is exponentially distributed with
parameter $2\sigma^{-2}(C_j-\sum_i A_{ji} \rho_i)$ 
for each $j \in \J$. 
The Brownian model
line  sizes and delays are again given by 
equations~(\ref{delay1}) and (\ref{delay2}) respectively, each
with stationary distributions given by  a linear combination
of independent exponential random variables, one for each
section of road. 

A key feature of the linear network, and its generalization
to tree networks, is that all traffic is bound for the same destination.
In our application to a road network this ensures that all
traffic in a  line at a given entry point is on the same route.
If traffic on different routes shared a single  line it
would not be possible to align the delay incurred by
traffic so precisely with the sum of dual variables 
for the resources to be used.\footnote{The tree 
topology of Figure~\ref{fig:roadtree} ensures that the queueing delays 
in the proportionally fair Brownian network model 
are partially ordered.
A technical consequence
is that a wide class of fair
\index{network!capacity|(}capacity allocations, 
the $\alpha$-fair allocations, share the same workload cone:
in the notation
\index{Kang, W. N.|pagenote}\index{Lee, N. H.|pagenote}of~\cite{KKLW}, the
\index{workload!cone|)}cone $\mathcal W_\alpha$ 
does not depend upon $\alpha$.}

\section{Route choices} \label{rc}

\begin{figure}[!t]
\begin{center}
\psfrag{W1}[t][t]{$w_1$}
\psfrag{W2}[t][t]{$w_2$}
\psfrag{0}[t][t]{$0$}
\psfrag{User}[t][t]{users}
\psfrag{long}[t][t]{long, uncongested link}
\psfrag{short}[t][t]{short,}
\psfrag{r1}[t][t]{$\rho_1$}
\psfrag{r2}[t][t]{$\rho_2$}
\psfrag{r3}[t][t]{$\rho_3$}
\psfrag{r4}[t][t]{$\rho_4$}

\psfrag{C1}[t][t]{$ C_1$}
\psfrag{C2}[t][t]{$ C_2$}
\psfrag{C3}[t][t]{$ C_3$}
\psfrag{C4}[t][t]{$ C_4$}

\psfrag{x1}[t][t]{$m_1$}
\psfrag{x2}[t][t]{$m_2$}
\psfrag{x3}[t][t]{$m_3$}
\psfrag{x4}[t][t]{$m_4$}

\centerline\leavevmode\epsfig{file=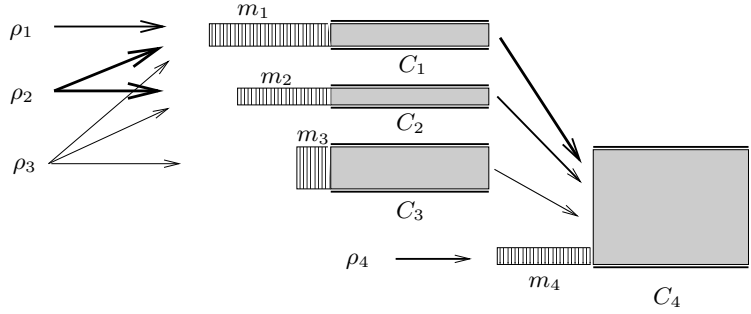,width=3.8in}

\end{center}
\caption{
Three parallel roads lead to a fourth road and hence to a 
common destination.
Lines of size $m_1$, $m_2$, $m_3$, $m_4$ are held on the slip roads leading
to these roads.  
There are four sources of traffic:  sources $2$ and $3$ 
may choose their first road, with choices as shown.
}
 \label{fig:roadparallel}
\end{figure}

Next consider the road network illustrated in
Figure~\ref{fig:roadparallel}. Three parallel roads lead
into a fourth road and hence to a common destination.
Access to each of these roads
is metered, so that their
respective capacities $C_1$, $C_2$, $C_3$, $C_4$
are not overloaded, and $C_1+ C_2+ C_3 < C_4$. 
There are four sources of traffic with
respective loads $\rho_1$, $\rho_2$, $\rho_3$, $\rho_4$: the first source
has access to road $1$ alone, on its way to road $4$; 
the second source has access to both
roads $1$ and $2$; and the third source can access all three 
of the parallel roads.
We assume that traffic arriving with access to more than
one road  distributes itself in an attempt to minimize its
queueing delay, an assumption whose implications 
we shall explore. 

We could view sources of traffic as arising in different
geographical regions, with different possibilities for easy 
access to the motorway network and with real time information
on delays. Or we could imagine a priority access
\index{queue!discipline}discipline
where some traffic, for example high occupancy vehicles,
has a larger set of lines to choose from. 

Given the line sizes
$M(t)=m$, we suppose the metered rates $\Lambda(m)$ are
chosen to be proportionally
\index{fairness|)}fair: that is, the
\index{network!capacity allocation policy}capacity allocation policy
$\Lambda(\cdot)$ solves the
\index{convex optimization}optimization problem~(\ref{pf1})--(\ref{pf4}).
For this network 
\begin{equation*}
 A = \left( \begin{array}{cccc}
1 &  0  & 0 & 0 \\
0  &  1 & 0  & 0\\
0 & 0 & 1 & 0 \\
1  &   1       & 1 & 1
\end{array} \right), 
\qquad
C=
 \left( \begin{array}{c}
C_1 \\
{C}_2 \\
C_3 \\
C_4
\end{array} \right)
\end{equation*}
and so, from relations~(\ref{cs1})--(\ref{cs2}), 
\begin{equation*}
\Lambda_i(m) =  \frac{m_i}{q_i+ q_4 },  \, \ i=1,2,3, \qquad
\Lambda_4(m) =  \frac{m_4}{q_4}. 
\end{equation*}
We assume the ramp metering policy has no knowledge of 
the routing choices available to arriving traffic,
but is simply a function of the observed line  
sizes $m$, the topology matrix $A$ and the
\index{network!capacity|)}capacity
vector $C$. 

How might arriving traffic choose between lines? Well,
traffic that arrives when the line sizes are $m$ and
the metered rates are $\Lambda(m)$ might reasonably
consider the ratios $m_i/\Lambda_i(m)$
in order to choose which  line to join, since these
ratios give the time it would take 
to process the work currently
in line $i$ at the current metered rate for line  $i$, for $i=1$, 2, 3.
But these ratios are just $q_i+ q_4$ for $i=1$, 2, 3.
Given the choices available to the three sources, we would
expect 
exercise of these choices to ensure 
that $q_1 \geq q_2 \geq q_3$, or equivalently
that the delays through lines 1, 2, 3 are weakly decreasing. 




Because traffic from sources $2$ and $3$ has the ability to 
make route choices, condition~(\ref{capcon}) is sufficient,
but no longer
necessary, for stability. The stability condition for the
network of Figure~\ref{fig:roadparallel} is 
\begin{equation} \label{enl}
\sum_{i=1}^j \rho_i   < \sum_{i=1}^j C_i, \, \ j=1,2,3, \qquad
\sum_{i=1}^4 \rho_i  < C_4,
\end{equation}
and is thus of the form~(\ref{capcon}), but with $A$ and
$C$ replaced by $\bar{A}$ and $\bar{C}$ respectively, where 
\[ \bar{A} = \left( \begin{array}{cccc}
1 &  0  & 0 & 0 \\
1  &  1 & 0 & 0\\
1  &   1  & 1 & 0\\
1  &   1  & 1 & 1
\end{array} \right), 
\qquad
\bar{C}= 
 \left( \begin{array}{c}
\bar{C}_1 \\
\bar{C}_2 \\
\bar{C}_3 \\
\bar{C}_4
\end{array} \right) 
=
 \left( \begin{array}{c}
C_1 \\
 C_1 +  C_2 \\
C_1 +  C_2 +  C_3  \\
C_4 
\end{array} \right). 
\]
The forms $\bar{A}$, $\bar{C}$ capture the concept of 
four
\index{virtual resource|(}\emph{virtual resources} of capacities
$\bar{C}_j$, $j=1$, 2, 3, 4. 
Given the line  
sizes $m = (m_i, i=1,2,3,4)$, the workloads 
$w = (w_i, i=1,2,3,4)$ for the
four virtual resources are $w=\bar{A} m$.

For $i=1$, 2, 3, 4, 
we model  the cumulative inflow of work  from
source $i$ over the interval $(0,t]$ as a
Brownian motion $\tE_i=(\tE_i(t), t\geq 0)$ starting from
the origin with drift  $\rho_i$ and variance parameter $\rho_i \sigma^2$
that  can be written in the form~(\ref{eT}),
where $(\tZ_i(t), t \geq 0)$, $i=1$, 2, 3, 4, are independent
standard Brownian motions. Let $\tE = (\tE_i, i=1,2,3,4)$
and let   $\tX = (\tX_j, j \in \J)$ 
be defined by
\begin{equation*} 
\tX(t) =  \bar{A}  \tE(t) - \bar{C} t, \quad t\geq 0. 
\end{equation*}
Then $\tX$ is
a  four-dimensional Brownian motion
starting from
the origin with drift $  \bar{A}\rho-\bar{C}$ and covariance
matrix $\Gamma= \sigma^2 \bar{A}  \,  [\rho] \bar{A}'$.
We assume the stability 
condition~(\ref{enl}) is satisfied, so that all components
of the drift are strictly
 negative.  Let $\mathcal W, \mathcal W^j$ be
defined by~(\ref{wc}), (\ref{wcj}) respectively, with 
$A$ replaced by $\bar{A}$. 

Our Brownian network model for $\bar AM$ 
is then the process $\tW$
defined by properties (i)--(iv) of Section~\ref{bnm}
with $\tW$ in place of $\tilde W$ and $\bar C-\bar A\rho$ in
place of $\theta$. 

Define a process of dual variables for the
virtual resources:
 $\tQ = (A[\rho]A')^{-1} \tW$.
Since $\tW $ is our model for $\bar AM$,
our Brownian model for the line sizes is
given by
\begin{equation*}
\tM= \bar A^{-1} \tW =   [\rho]  \bar{A}' \tQ. 
\end{equation*}
Our Brownian model for the delays at each line is given by 
\begin{equation*}
\tD =  \bar{A}' \tQ, 
\end{equation*}
which from the form of $\bar{A}$ becomes
\begin{equation*}
\tD_i = \sum_{j=i}^{4} \tQ_j ,   \qquad i=1,2,3 ,4.
\end{equation*}
Thus at any time $t$ when $\tQ_1(t)>0$, $\tQ_2(t)>0$ then
$\tD_1(t) >\tD_2(t)>\tD_3(t)$, and the incentives for arriving traffic are
such that traffic from source $i$ joins line $i$.
However if $\tQ_1(t)=0$, and so $\tD_1(t) = \tD_2(t)$,
then arriving traffic from stream $2$ may choose to enter
line $1$, and thus contribute to increments of
the workload for virtual resource $1$, 
whilst still contributing to 
the workload for virtual resources 2, 3 and $4$. 
Our model corresponds
to the assumption that no more traffic does this than is
necessary to keep $\tQ_1$ non-negative, or equivalently
to keep $\tD_1 \geq \tD_2$.
Similarly if $\tQ_2(t)=0$
then $\tD_2(t) = \tD_3(t)$, and
arriving traffic from stream $3$ may choose to enter
line $2$, and thus  contribute to increments of
the workload for virtual resource $2$, 
whilst still contributing to
the workload for virtual
resources $3$ and $4$; we suppose
just sufficient traffic does this 
to keep $\tQ_2$ non-negative, or equivalently
to keep
$\tD_2 \geq \tD_3$. Finally if
$\tQ_3(t)=0$ or $\tQ_4(t)=0$ then (real)
resource $3$ or $4$ respectively may not be fully
utilized, as in earlier examples, and our model
 corresponds
to the assumption that there is no more loss of
utilization at (real) resources $3$ and $4$
than is necessary to prevent $\tW$ from
leaving $\mathcal W$.

If $\tW^s$ is a random variable with the
stationary distribution of $\tW$, then
the components of $\tQ^s = (A[\rho]A')^{-1} \tW^s$
are independent and for $j=1$, 2, 3, 4, $\tQ_j^s$ 
is exponentially distributed with
parameter $\zeta_j$ where
\begin{eqnarray*}
\zeta_1 &=&\frac{2}{\sigma^2}(  C_1 - \rho_1 ),
\\
\zeta_2 &=& \frac{2}{\sigma^2}( C_1 +C_2 - \rho_1 -\rho_2),
\\
\zeta_3  &=&\frac{2}{\sigma^2}(  C_1 +C_2 +C_3 - \rho_1
-\rho_2 - \rho_3 ), \\
\zeta_4  &=& \frac{2}{\sigma^2}( C_4 -  \rho_1
-\rho_2 - \rho_3 - \rho_4 ). 
\end{eqnarray*}
Under the Brownian
network model, the stationary distribution for line sizes and for delays
at each line are given by the distributions
of $\tM^s$ and $\tD^s$, respectively, where
\begin{equation*}
\tM^s_i =   \rho_i \left(\sum_{j=i}^4 \tQ^s_j \right) , 
\quad
 \tD_i^s= \sum_{j=i}^4 \tQ^s_j  
\quad i=1,2,3,4. 
\end{equation*}




The Brownian network model thus corresponds
to natural assumptions about how arriving
traffic from different sources would choose their routes. 
The results on the stationary distribution for 
the network are intriguing. 
The ramp metering policy
has no knowledge of the routing
choices available to arriving traffic, and hence of the
enlarged stability region~(\ref{enl}). 
Nevertheless, under the Brownian model,
 the interaction of
the ramp metering policy with the routing
choices available to arriving traffic has a performance
described in terms of 
dual random variables, one for each of
the virtual resources of the enlarged stability region;
when a driver makes a route choice,
the delay facing a driver on a route is 
a  sum of dual random variables, one for each of the 
\index{virtual resource|)}virtual resources used by that route; and under
their stationary distribution, the dual random variables
are independent and exponentially distributed.\index{network|)}

\section{Concluding remarks}

The design of
\index{ramp metering|)}ramp metering strategies cannot assume that arriving
traffic flows are exogenous, since in general drivers' behaviour will be
responsive to the delays incurred or expected.  In this paper we have
presented a preliminary exploration of an approach to the design of ramp
metering flow rates informed by earlier work on
\index{Internet}Internet congestion control. A feature of this approach is
that it may prove possible to integrate ideas of
\index{fairness}fairness of a control policy
with overall system optimization.

There remain many areas for further investigation.  In particular, we have
seen intriguing examples, in the context of a single queue and of Internet
congestion control, of remarkably good approximations produced for the
stationary distributions of queue length and workload by use of the direct
Brownian modelling approach. Furthermore, in the context of a controlled
motorway, where a detailed model for arriving traffic is not easily available,
use of a direct Brownian model has enabled us to develop an approach to the
design and performance of ramp metering and in the context of that model to
obtain insights into the interaction of ramp metering with route
choices.\index{road traffic|)} Nevertheless, we expect that the use of direct
\index{Brown, R.!Brownian network|)}Brownian network models will not always
produce good results.  Indeed, it is
possible that such models may be suitable only when the scaled
\index{workload|)}workload
process can be approximated in heavy traffic by a reflecting
\index{Brown, R.!Brownian motion|)}Brownian motion 
that has a product-form stationary distribution. 
We believe that understanding when the direct method
is a good modelling approach  and when it is not,  
and obtaining a rigorous  understanding of the
reasons for this, is  an interesting
topic worthy of  further research.\index{heavy traffic|)}


\end{document}